\documentclass[%
 aip,
cp,  
 amsmath,amssymb,
 reprint,%
]{revtex4-2}

\usepackage{graphicx}
\usepackage{epstopdf}
\usepackage{dcolumn}
\usepackage{bm}
\usepackage[utf8]{inputenc}
\usepackage[T1]{fontenc}
\usepackage{mathptmx}

\begin{document}

\title{Application of the Adams-Bashfort-Mowlton Method to the Numerical Study of Linear Fractional Oscillators Models}

\author{Parovik Roman} 
 \email[Corresponding author: ]{romanparovik@gmail.com}
\affiliation{Institute of Cosmophysical Research and Radio Wave Propagation, Far Eastern Branch of the Russian Academy of Sciences
  (use complete addresses, including country name or code).
}

\date{\today} 

\begin{abstract}
The paper presents a numerical analysis of the class of mathematical models of linear fractional oscillators, which is the Cauchy problem for a differential equation with derivatives of fractional orders in the sense of Gerasimov-Caputo. A method based on an explicit nonlocal finite-difference scheme (ENFDS) and the Adams-Bashfort-Moulton (ABM) method are considered as tools for numerical analysis. An analysis of the errors of the methods is carried out, it is shown that the ABM method is more accurate and converges faster to an exact solution than the ENFDS method.
\end{abstract}

\maketitle

\section{Introduction}

Models of oscillatory systems (oscillators) are used in various fields of knowledge from mechanics to economics and biology \cite{Adronov_2013,Danylchuk_2019,Li_2018}. From the point of view of mathematics, these models are usually described using ordinary differential equations of the second order and the corresponding initial conditions, i.e. the Cauchy problem is posed \cite{Adronov_2013}. It should be noted that such mathematical models cannot take into account the properties of the environment, for example, heredity (heredity) or memory effects. This effect is characterized by the fact that the oscillating medium “remembers” the impact on it for a long time.

For the first time, the model of the hereditary oscillator was presented in his work by the Italian mathematician V. Volterra \cite{Volterra_1912}. He proposed to take into account heredity in the linear oscillator model using an integro-differential equation with a difference kernel, which he later called the function of heredity or memory. It should be noted that V. Volterra derived the total energy conservation law for this generalized oscillator, in the formula of which an additional term appeared, which is responsible for the dissipation of its energy. This important result was confirmed in subsequent works on this topic.

If we choose a power-law memory function, then we can, using the mathematical apparatus of fractional calculus, go to other model equations that contain derivatives of fractional orders. In this case, the orders of fractional derivatives, as shown by the results of \cite{Pskhu_2018,Parovik_2019}, will be responsible for the intensity of energy dissipation and are related to the Q-factor of the oscillator. Oscillators with such a description are usually called fractional.

Research methods for mathematical models of fractional oscillators can be divided into exact and numerical. The exact methods, for example, include integral transformations \cite{Khalouta_2019} or decomposition methods \cite{Momani_2007}, and numerical methods include the theory of finite-difference schemes \cite{Diethelm_2005}, variational-iterative methods \cite{Momani_2007}.

In this paper, we will carry out a numerical analysis of mathematical models of fractional linear oscillators (FLO) using elements of the theory of finite-difference schemes. As methods of numerical analysis, we will choose a method based on an explicit nonlocal finite-difference scheme studied in the author's work \cite{Parovik_2018} and the fractional ABM method, which was investigated in \cite{Garrappa_2018,Yang_2006,Diethelm_2002}. Let us analyze the errors of the methods using Runge's rule.

\section{Some reduction from the theory of fractional calculus}

Here we will consider the main definitions from the theory of fractional calculus, in more detail its aspects can be studied in the books \cite{Kilbas_2006,Oldham_1974,Miller_1993}.

\textbf{Definition 1.} Fractional Riemann-Liouville integral of order $\alpha$:
	\begin{equation}
		I_{0t}^\alpha x\left( \tau  \right) = \frac{1}{{\Gamma \left( \alpha 
				\right)}}\int\limits_0^t {\frac{{x\left( \tau  \right)d\tau }}{{{{\left( {t -
								\tau } \right)}^{1 - \alpha }}}}} ,\alpha  > 0,t > 0,
		\label{Parovikeq1}
\end{equation}	
here $\Gamma\left(.\right)$~ is the gamma function.

The operator (\ref{Parovikeq1}) has the following properties:
\[
I_{0t}^0x\left( \tau  \right) = x\left( t \right), I_{0t}^\alpha I_0^\beta
x\left( \tau  \right) = I_{0t}^{\alpha  + \beta }x\left( \tau  \right), I_{0t}^\alpha I_0^\beta x\left( \tau  \right) = I_0^\beta I_{0t}^\alpha x\left(\tau  \right).
\]	

\textbf{Definition 2.} The fractional Gerasimov-Caputo derivative of order $\alpha$ has the form:
	\begin{equation}
		\partial _{0t}^\alpha x\left( \tau  \right) = \left\{ \begin{array}{l}
			\dfrac{1}{{\Gamma \left( {m - \alpha } \right)}}\int\limits_0^t
			{\dfrac{{{x^{(m)}}\left( \tau  \right)d\tau }}{{{{\left( {t - \tau }
								\right)}^{\alpha  + 1 - m}}}},0 \le m - 1 < \alpha  < m,} \\
			\dfrac{{{d^m}x\left( t \right)}}{{d{t^m}}},m \in N.
		\end{array} \right.
		\label{Parovikeq2}
	\end{equation}		

The operator (\ref{Parovikeq2}) has the following properties\cite{Kilbas_2006}:
\[
\partial _{0t}^\alpha I_{0t}^\alpha x\left( \tau  \right) = x\left( t
\right),I_{0t}^\alpha \partial _{0t}^\alpha x\left( \tau  \right) -
\sum\limits_{k = 0}^{i - 1} {\frac{{{x^{\left( k \right)}}\left( 0
			\right){t^k}}}{{k!}}} ,t > 0.
\]

\section{Statement of the problem}

Consider the following Cauchy problem:

\begin{equation}
	\partial _{0t}^\beta x\left( \tau  \right) + \lambda \partial _{0t}^\gamma x\left( \tau  \right) + \omega \left( t
	\right)x\left( t \right) = f\left( t \right), x\left( 0 \right) = {\alpha_1},\dot x\left( 0 \right) = {\alpha _2},
	\label{Parovikeq3}
\end{equation}
where $ x\left(t \right) \in {C^2}\left[{0, T} \right] $ is the solution (displacement) function; $ t \in \left [{0, T} \right] $ is the coordinate, responsible for time, $ T> 0 $ is the constant, simulation time; $ \omega \left(t \right) $ is the function responsible for the frequency of free oscillations and determines the type of FLO; $ \lambda $ is the coefficient of friction; $ f \left(t \right)$ is the function responsible for external influence; $ {\alpha _1}, {\alpha _2} $ is the given constants defining the initial condition, fractional differentiation operators are understood in the Gerasimov-Caputo sense (\ref{Parovikeq2}) of orders $ 1 <\beta <2 $ and $ 0 <\gamma <1 $.

\textbf{Remark 1.} The Cauchy problem (\ref{Parovikeq3}) describes a wide class of FLO and in the case of $ \beta = 2 $ and $ \gamma = 1 $ becomes the class of ordinary linear oscillators \cite{Adronov_2013}.

\textbf{Definition 3.} The Cauchy problem (\ref{Parovikeq3}) will be called the fractional linear oscillator model (FLO).

\section{Methods of solution}


Consider two numerical methods for solving the Cauchy problem (\ref{Parovikeq3}): ENFDS and ABM method. The ABM method uses a combination of the explicit Adams-Bashforth method and the implicit Adams-Moulton method.
	
\textbf{Explicit non-local finite-difference scheme.} Let $x\left(t\right) \in C^3\left[0,T\right]$ to achieve the required smoothness. Divide the time interval $ [0, T] $ into $ N $ equal parts with a constant step $\tau = {T \mathord{\left/{\vphantom {TN}} \right. \kern- \nulldelimiterspace} N} $. The solution function $ x \left (t \right) $ will go to the grid function $ x \left({{t_k}}\right) = {x_k} $, $ k = 1, \ldots, \, N-1 $.

The Gerasimov-Caputo operators in the equation (\ref{Parovikeq3}) are approximated as follows. For the second operator from (\ref{Parovikeq3}):
\[
\partial _{0t}^\gamma x\left( \tau  \right) = \frac{1}{{\Gamma \left( {1 -
			\gamma } \right)}}\int\limits_0^t {\frac{{\dot x\left( \tau  \right)d\tau
	}}{{{{\left( {t - \tau } \right)}^\gamma }}}}  \approx \frac{1}{{\Gamma \left( {1
			- \gamma } \right)}}\sum\limits_{j = 0}^k {\int\limits_{j\tau }^{\left( {j + 1}
		\right)\tau } {\frac{{\dot x\left( \tau  \right)d\tau }}{{{{\left( {{t_{k + 1}} -
							\tau } \right)}^\gamma }}}} }  =
\]

\begin{equation}
	=\frac{1}{{\Gamma \left( {1 - \gamma } \right)}}\sum\limits_{j = 0}^k
	{\frac{{x\left( {{t_{j + 1}}} \right) - x\left( {{t_j}} \right)}}{\tau
		}\int\limits_{j\tau }^{\left( {j + 1} \right)\tau } {\frac{{d\tau }}{{{{\left(
							{{t_{k + 1}} - \tau } \right)}^\gamma }}}} }  = 
	\label{Parovikeq4}
\end{equation}
\[
\frac{1}{{\Gamma \left( {1 -
			\gamma } \right)}}\sum\limits_{j = 0}^k {\frac{{x\left( {{t_{j + 1}}} \right) -
			x\left( {{t_j}} \right)}}{\tau }\int\limits_{\left( {k - j} \right)\tau }^{\left(
		{k - j + 1} \right)\tau } {\frac{{d\eta }}{{{\eta ^\gamma }}}} }  =\frac{1}{{\Gamma \left( {1 - \gamma } \right)}}\sum\limits_{j = 0}^k
{\frac{{x\left( {{t_{k - j + 1}}} \right) - x\left( {{t_{k - j}}} \right)}}{\tau
	}\int\limits_{j\tau }^{\left( {j + 1} \right)\tau } {\frac{{d\eta }}{{{\eta
					^\gamma }}}} }  = 
\]
\[				
= \frac{{{\tau ^{ - \gamma }}}}{{\Gamma \left( {2 - \gamma }
		\right)}}\sum\limits_{j = 0}^k {\left( {{{\left( {j + 1} \right)}^{2 - \beta }} -
		{j^{2 - \beta }}} \right)\left( {x\left( {{t_{k - j + 1}}} \right) - x\left(
		{{t_{k - j}}} \right)} \right)} .
\]

Therefore, we have according to (\ref{Parovikeq4}):
\begin{equation}
	\partial _{0t}^\gamma x\left( \tau  \right) \approx B\sum\limits_{j= 0}^{k - 1} {{b_j}\left( {{x_{k - j + 1}} - {x_{k - j - 1}}} \right)}, {b_j} = {\left( {j + 1} \right)^{1 - \gamma }} - {j^{1 - \gamma }}, B = \frac{{\lambda{\tau ^{ - \gamma }}}}{{\Gamma \left( {2 - \gamma } \right)}},
	\label{Parovikeq5}
\end{equation}

Similarly, for the first operator from (\ref{Parovikeq3}) we will have:

\begin{equation}
	\partial _{0t}^\beta x\left( \eta  \right) \approx A\sum\limits_{j = 0}^{k - 1}
	{{a_j}\left( {{x_{k - j + 1}} - 2{x_{k - j}} + {x_{k - j - 1}}} \right)} ,{a_j} =
	{\left( {j + 1} \right)^{2 - \beta }} - {j^{2 - \beta }},
	\label{Parovikeq6}
\end{equation}

Substituting the approximations (\ref{Parovikeq5}) and (\ref {Parovikeq6}) in the original equation (\ref{Parovikeq3}), we arrive at the following discrete Cauchy problem:
\begin{equation}
	A\sum\limits_{j = 0}^{k - 1} {{a_j}\left({{x_{k - j + 1}} - 2{x_{k - j}} + {x_{k - j - 1}}} \right)}  + B\sum\limits_{j =
		0}^{k - 1} {{b_j}\left( {{x_{k - j + 1}} - {x_{k - j}}} \right)}  + {\omega_k}{x_k} = {f_k},
	\label{Parovikeq7}           
\end{equation}
\[
{x_0} = {\alpha _1},{x_1} = {\alpha _2} + \tau {\alpha _1}.
\]

For the discrete Cauchy problem (\ref{Parovikeq7}), the following lemma holds.

\textbf{Lemma 1.} The coefficients of the discrete Cauchy problem (\ref{Parovikeq7}) have the following properties:
	\begin{equation}
		1) \sum\limits_{j = 0}^{k - 1} {{a_j} = {k^{2 - \beta }}},\, \sum\limits_{j = 0}^{k-1} {{b_j} = {k^{1 - \gamma }}},\, 
		2)\, 1 = a_0>a_1>\ldots>0, 1 = b_0>b_1>\ldots>0,\, 3)\, A\geq 0, B \geq 0
		\label{Parovikeq8}
	\end{equation}

\textbf{Proof.}	The first property follows from the definition:
\[
\sum\limits_{j = 0}^{k - 1} {{a_j} = \sum\limits_{j = 0}^{k - 1} {\left[
		{{{\left( {j + 1} \right)}^{2 - \beta }} - {j^{2 - \beta }}} \right]} }  = 1 - 0
+ {2^{2 - \beta }} - 1 + {3^{2 - \beta }} - {2^{2 - \beta }} + ... +
\]

\[
+ {\left( {k - 1} \right)^{2 - \beta }} + {k^{2 - \beta }} - {\left( {k - 1}
	\right)^{2 - \beta }} = {k^{2 - \beta }}.
\]

\[
\sum\limits_{j = 0}^{k - 1} {{b_j} = \sum\limits_{j = 0}^{k - 1} {\left[
		{{{\left( {j + 1} \right)}^{1 - \gamma }} - {j^{1 - \gamma }}} \right]} }  = 1 -
0 + {2^{1 - \gamma }} - 1 + {3^{1 - \gamma }} - {2^{1 - \gamma }} + ... +
\]

\[
+ {\left( {k - 1} \right)^{1 - \gamma }} + {k^{1 - \gamma }} - {\left( {k - 1}
	\right)^{1 - \gamma }} = {k^{1 - \gamma }}.
\]	

We prove the second condition as follows. Let's introduce the following functions: $\phi\left(z\right)=(z+1)^{2-\beta}-z^{2-\beta}$, $\eta\left(z\right)=(z+1)^{1-\gamma}-z^{1-\gamma}$ where $z > 0$. These functions are monotonically decreasing. Indeed, let us take the derivatives with respect to $z$ of these functions. We get:
$$ \phi'\left(z\right)=\left(2-\beta\right)(z+1)^{1-\beta}-z^{1-\beta},\, \eta'\left(z\right)=\left(1-\gamma\right)(z+1)^{1-\gamma}-z^{1-\gamma}$$

The third property follows from the definition of the gamma function. It is known that the gamma function $ \Gamma\left(z\right)> 0 $ for $ \forall z $. Since $ \tau> 0, \lambda> 0 $, then we come to the conclusion that $ A \geq 0 $, $ B \geq 0 $. \hfill $\square$

Let us now investigate the order of approximation of fractional operators $\partial_{0t}^\beta x\left(\tau\right)$ and $\partial_{0t}^\gamma x\left(\tau\right)$. Let $\overline{\partial}_{0t}^\beta x\left(\tau\right)$ and $\overline{\partial}_{0t}^\gamma x\left(\tau\right)$ is the approximation operators. Then we have the following lemma.

\textbf{Lemma 2.} Approximations $\overline{\partial}_{0t}^\beta x\left(\tau\right)$ and $\overline{\partial}_{0t}^\gamma x\left(\tau\right)$  Gerasimov-Caputo operators (\ref{Parovikeq2}) $\partial_{0t}^\beta x\left(\tau\right)$ and $\partial_{0t}^\gamma x\left(\tau\right)$ satisfy the following estimates:
	\begin{equation}
		\partial_{0t}^\beta x\left(\tau\right) =\overline{\partial}_{0t}^\beta x\left(\tau\right)+ O\left(\tau^{4-\beta}\right), \partial_{0t}^\gamma x\left(\tau\right) =\overline{\partial}_{0t}^\gamma x\left(\tau\right)+ O\left(\tau^{2-\gamma}\right),
	\label{Parovikeq9}
	\end{equation}
where $O\left(\cdot\right)$ is the Landau symbol.

\textbf{Proof}. Using the first property (\ref {Parovikeq8}) of Lemma 1 and the relations (\ref{Parovikeq5}) and (\ref{Parovikeq6}), we obtain
\[
\overline{\partial}_{0t}^\beta x\left(\tau\right)=\frac{\tau^{2-\beta}}{\Gamma\left(3-\beta\right)}\sum\limits_{j=0}^{k-1}a_j\left[\ddot{x}\left(t-j\tau\right)+O\left(\tau^2\right)\right]=\frac{\tau^{2-\beta}}{\Gamma\left(3-\beta\right)}\sum\limits_{j=0}^{k-1}a_j\ddot{x}\left(t-j\tau\right)+\frac{k^{2-\beta}O\left(\tau^{4-\beta}\right)}{\Gamma\left(3-\beta\right)}=
\]
\[
=\frac{\tau^{2-\beta}}{\Gamma\left(3-\beta\right)}\sum\limits_{j = 0}^{k - 1} {{a_j}\ddot x\left( {t - j\tau } \right) + O\left( {{\tau^{4-\beta}}} \right)}.
\]
\[
\partial _{0t}^\beta x\left( \tau  \right) = \frac{{{\tau ^{2 - \beta
}}}}{{\Gamma \left( {3 - \beta } \right)}}\sum\limits_{j = 0}^{k - 1}
{{a_j}\left[ {\ddot x\left( {t - {\xi _j}} \right)} \right]} ,{\xi _j} \in \left[
{j\tau ,\left( {j + 1} \right)\tau } \right].
\]
\[
|\overline{\partial}_{0t}^\beta x\left(\tau\right)-\partial_{0t}^\beta x\left(\tau\right)|=\left|\frac{\tau^{2-\beta}}{\Gamma\left(3-\beta\right)}\sum\limits_{j = 0}^{k-1}a_j\left[\ddot{x}\left(t-j\tau\right)-\ddot{x}\left(t-\xi_j\right)\right]+O\left(\tau^{4-\beta}\right)\right|=
\]
\[
=\left|\frac{\tau^{2-\beta}}{\Gamma\left(3-\beta\right)}\sum\limits_{j = 0}^{k-1}a_j O\left(\tau^2\right)+O\left(\tau^{4-\beta}\right)\right|=\left|\frac{k^{2-\beta}O\left(\tau^{4-\beta}\right)}{\Gamma\left(3-\beta\right)}+O\left(\tau^{4-\beta}\right)\right|=O\left(\tau^{4-\beta}\right). 
\]
Likewise, we can show the second estimate (\ref {Parovikeq9}). The lemma is proved. \hfill $\square$

\textbf{Remark 2.} In the case when in relations (9) $ \beta = 2 $ and $\gamma = 1 $, then we obtain approximations of the derivatives of the first and second orders with the corresponding first and second orders of approximation.

\textbf{Remark 3.} Note that at the internal nodes of the computational grid, taking into account Lemma 2, the explicit nonlocal finite-difference scheme (7) approximates the differential Cauchy problem with order $ 2- \gamma $, however, the general order of approximation of scheme (7) is the first due to the approximation of the initial conditions. 

We rewrite the discrete Cauchy problem (\ref{Parovikeq7}) as follows:
\begin{equation}
\left\{\begin{array}{l}
{x_0} = {\alpha _1},\, j =0,\\
{x_1} = {\alpha _1} + \tau {\alpha _2},\, j = 1,\\
{x_{k + 1}} = \dfrac{{{f_k} + \left({2A - B - {\omega _k}} \right){x_k} - A{x_{k- 1}}}}{{A + B}} - \\ - \dfrac{A}{{A + B}}\sum\limits_{j = 1}^{k - 1} {{a_j}\left( {{x_{k - j + 1}}-2{x_{k - j}} + {x_{k - j - 1}}} \right)}- \\ - \dfrac{B}{{A + B}}\sum\limits_{j =1}^{k - 1} {{b_j}\left( {{x_{k - j + 1}} - {x_{k - j}}} \right)},\,j = 2,\ldots,\,N -1.
\end{array} \right.
\label{Parovikeq11}
\end{equation}
or in matrix form:
\begin{equation}
{X_{k + 1}} = M{X_k} + {F_k},
\label{Parovikeq12}
\end{equation}
\[
{X_{k + 1}} = {\left( {{x_1},{x_2},\ldots,{x_k}} \right)^T},\,{X_k} = {\left({{x_0},{x_1},\ldots,{x_{k - 1}}} \right)^T},
\] 
\[
{F_k} = {\left( {{f_0},\frac{{{f_1}}}{{A + {\lambda _1}B}},\ldots,\frac{{{f_{k -1}}}}{{A + {\lambda _{k - 1}}B}}} \right)^T},\,{f_0} = \tau {\alpha_2},
\]
where the Hessenberg matrix $ M $ in (\ref{Parovikeq12}) has the form:
\begin{equation}
{m_{ij}} = \left\{ \begin{array}{l}
0,\,j \ge i + 1,\\
\dfrac{{A\left( {2 - {a_1}} \right) - B{\lambda _{i - 1}}{b_1}}}{{A + {\lambda_{i - 1}}B}},\,j = i = 3,\ldots,\,N - 1,\\
\dfrac{{ - A\left( {{a_{i - j + 1}} - 2{a_{i - j}} + {a_{i - j - 1}}} \right)-B{\lambda _{i - 1}}\left( {{b_{i - j + 1}} - {b_{i - j - 1}}} \right)}}{{A +{\lambda _{i - 1}}B}},\,j \le i - 1, \end{array} \right.
\label{Parovikeq13}
\end{equation}
\[
{m_{1,1}} = 1,\,{m_{2,2}} = \dfrac{{2A}}{{A + {\lambda _1}B}},\,{m_{i,1}} =\dfrac{{B{\lambda _{i - 1}}{b_{i - 2}} - A{a_{i - 2}}}}{{A + {\lambda _{i-1}}B}},\,i = 2,\ldots,\,N - 1,
\]
\[
{m_{i,2}} = \dfrac{{A\left( {2{a_{i - 2}} - {a_{i - 3}}} \right) + {\lambda _{i-1}}B{b_{i - 3}}}}{{A + {\lambda _{i - 1}}B}},\,i = 3,\ldots,\,N - 1.
\]

The following theorems are true.

\textbf{Theorem 1 \cite{Parovik_2018}}.
A nonlocal explicit finite-difference scheme (\ref{Parovikeq11}) converges with the first order if the following condition is satisfied:
	\begin{equation}
	\tau  \le {\tau _0} = \min \left( {1,{{\left( {\dfrac{{\Gamma \left( {2 - \gamma} \right)}}{{\lambda \Gamma \left( {3 - \beta } \right)}}}\right)}^{\dfrac{1}{{\beta  - \gamma }}}}} \right).
	\label{Parovikeq14}
	\end{equation}

Let $ {X_k}, {Y_k} $ be two different solutions of the matrix equation (\ref {Parovikeq12}) with the initial conditions $ {X_0}, {Y_0} $. Then the scheme stability theorem is valid.

\textbf{Theorem 2 \cite{Parovik_2018}}.
	A nonlocal explicit finite-difference scheme (\ref {Parovikeq11}) is conditionally stable if the condition (\ref {Parovikeq14}) is satisfied and the estimate $ \left | {{Y_k} - {X_k}} \right | \le C \left | {{Y_0} - {X_0}} \right | $ for any $ k $, where $ C> 0 $ is a constant independent of the step $ \tau $.	
	
The proof of Theorems 1 and 2 is based on the results of Lemmas 1 and 2.

\textbf {Adams-Bashfort-Moulton method.} The ABM method is a type of numerical predictor-corrector method for solving differential equations. It has been studied and discussed in detail in the papers \cite{Garrappa_2018,Yang_2006,Diethelm_2002}. Let's generalize this method for solving the Cauchy problem (\ref {Parovikeq3}). Taking into account Definitions 1 and 2, as well as the corresponding properties of the operators of fractional integration and differentiation, we write it in the form of a system:
\begin{equation}
\left\{ \begin{array}{l}
\partial _{0t}^{{\sigma _1}}x\left( \tau  \right) = y\left( t \right),\\
\partial _{0t}^{{\sigma _2}}y\left( \tau  \right) = f\left( t \right) - \lambda
y\left( t \right) - \omega \left( t \right)x\left( t \right),\\
x\left( 0 \right) = {\alpha _1},y\left( 0 \right) = {\alpha _2}
\end{array} \right.
\label{Parovikeq15}
\end{equation}
where ${\sigma _1} = \gamma $, ${\sigma _2} = \beta  - \gamma $. 

Let $ 0 <{\sigma _2} \le 1 $, i.e. $ \left\{\beta \right\} <\left\{\gamma \right\} $, \, $ \left\{\cdot \right\} $ is the fractional part of the number.

\textbf{Remark 4.} Note that if the system (\ref {Parovikeq15}) $ \left\{\beta\right\}> \left\{\gamma \right \} $, i.e. $ 1 <{\sigma_2} <2 $, then we can always transform it by adding one more equation.

On a uniform grid, we introduce the functions $ x_{n + 1}^p, \, y_{n + 1}^p $, $n = 0, \ldots, \, N-1 $, which will be determined by the Adams-Bashforth formula (predictor):
\begin{equation}
\left\{ \begin{array}{l}
x_{n + 1}^p = {x_0} + \dfrac{{{\tau ^{{\sigma _1}}}}}{{\Gamma \left( {{\sigma _1}
			+ 1} \right)}}\sum\limits_{j = 0}^n {\theta _{j,n + 1}^1{y_j},} \\
y_{n + 1}^p = {y_0} + \dfrac{{{\tau ^{{\sigma _2}}}}}{{\Gamma \left( {{\sigma _2}
			+ 1} \right)}}\sum\limits_{j = 0}^n {\theta _{j,n + 1}^2\left( { - \lambda {y_j}
		- {\omega _j}{x_j} + {f_j}} \right),} \\
\theta _{j,n + 1}^i = {\left( {n - j + 1} \right)^{{\sigma _i}}} - {\left( {n -
		j} \right)^{{\sigma _i}}},i = 1,2.
\end{array} \right.
\label{Parovikeq16}
\end{equation}

Then, using the Adams-Moulton formula for the corrector, we get:
\begin{equation}
\left\{ \begin{array}{l}
{x_{n + 1}} = {x_0} + \dfrac{{{\tau ^{{\sigma _1}}}}}{{\Gamma \left( {{\sigma _1}+ 2} \right)}}\left( {y_{n + 1}^p + \sum\limits_{j = 0}^n {\rho _{j,n +1}^1{y_j}} } \right),\\
{y_{n + 1}} = {y_0} + \dfrac{{{\tau ^{{\sigma _2}}}}}{{\Gamma \left( {{\sigma _2}+ 2} \right)}}\left( { - \lambda y_{n + 1}^p - {\omega _{n + 1}}x_{n + 1}^p + {f_{n + 1}} + \sum\limits_{j = 0}^n {\rho _{j,n + 1}^2\left( { - \lambda {y_j}-{\omega _j}{x_j} + {f_j}} \right)} } \right),\end{array} \right.
\label{Parovikeq17}
\end{equation}
where the weight factors in (\ref{Parovikeq17}) are determined by the formula:
\[
\rho _{j,n + 1}^i = \left\{ \begin{array}{l}
{n^{{\sigma _i} + 1}} - \left( {n - {\sigma _i}} \right){\left( {n + 1}\right)^{{\sigma _i}}},\,j = 0,\\ {\left( {n - j + 2} \right)^{{\sigma _i} + 1}} + {\left( {n - j}
	\right)^{{\sigma _i} + 1}} - 2{\left( {n - j + 1} \right)^{{\sigma _i} + 1}},\,1\le j \le n,\\ 1,\,j = n + 1,\\ i = 1,2. \end{array} \right.
\]

\textbf{Theorem 3 \cite{Yang_2006}.} If $\partial _{0t}^{{\sigma _i}}{x_i}\left( \tau  \right) \in {C^2}\left[ {0,T} \right],\left( {{x_1} = x\left( t \right),\,{x_2} = y\left(t\right),\,i = 1,2} \right)$, then
	\begin{equation}
	\mathop {\max }\limits_{1 \le j \le n} \left|{{x_i}\left( {{t_j}} \right) - {x_{i,j}}} \right| = O\left( {{h^{1 + \mathop
				{\min }\limits_i {\sigma _i}}}} \right).
	\label{Parovikeq18}
	\end{equation}
	
\section{Computational Accuracy Analysis}

Let us examine how the computational accuracy of the methods behaves. To do this, we will use the double recalculation method (Runge's rule) to estimate the error using the formula:
\begin{equation}
\xi  = \frac{{\mathop {\max}\limits_i \left| {{x_i} - {x_{2i}}} \right|}}{{{2^\mu } - 1}}, 
\label{Parovikeq19}
\end{equation}
where $ \mu $ is the order of approximation of the numerical method, $ {x_{2i}} $ is the numerical solution at the step $ {\tau\mathord{\left /{\vphantom {\tau 2}} \right. \kern- \nulldelimiterspace} 2} $. The computational accuracy of $ p $ is determined from the formula:
\begin{equation}
p = {\log _{\frac{{{\tau_1}}}{{{\tau _2}}}}}\frac{{{\xi _1}}}{{{\xi _2}}},
\label{Parovikeq20}                              
\end{equation}
$ {\tau _1}, \, {\tau _2} = {{{\tau _1}} \mathord {\left /{\vphantom {{{\tau _1}} 2}} \right. \kern- \nulldelimiterspace} 2} $ ~ --- grid steps, $ {\xi _1}, \, {\xi _2} $ ~ --- errors at step $ {\tau _1} $ and at step $ { \tau _2} $.

\textbf{Example 1}. Consider the Cauchy problem (\ref{Parovikeq3}) with homogeneous initial conditions, choosing the following functions and parameter values:
	$\omega \left( t \right) = \omega _0^\beta$,\, $f\left( t \right) = \omega _0^\beta {t^3} + \dfrac{{6{t^{3 - \beta }}}}{{\Gamma \left( {4 - \beta } \right)}} + \dfrac{{6\lambda {t^{3 - \gamma }}}}{{\Gamma \left( {4 - \gamma } \right)}},\,\lambda  = 0.1,\,\beta  = 1.8,\,\gamma  = 0.9,\, {\omega _0} = 1,\, t \in \left[ {0,1} \right]$.
\begin{equation}
	\partial _{0t}^{1.8}x\left( \tau  \right) + 0.1\partial _{0t}^{0.9}x\left( \tau 
	\right) + x\left( t \right) = {t^3} + \frac{{6{t^{1.2}}}}{{\Gamma \left( {2.2}
			\right)}} + \frac{{0.6{t^{2.1}}}}{{\Gamma \left( {3.1} \right)}},x\left( 0
	\right) = \dot x\left( 0 \right) = 0.
	\label{Parovikeq21}
\end{equation}
The solution to the Cauchy problem (\ref {Parovikeq21}), as you can see, is the function:
\begin{equation}
x\left( t \right)={t^3}. 
\label{Parovikeq22}                 
\end{equation}

\textbf{Remark 5.} The values of the parameters were chosen in such a way that the condition of Theorems 1 and 2 would be satisfied.

Due to the fact that the exact solution (\ref {Parovikeq22}) of the Cauchy problem (\ref {Parovikeq21}) is known, we can calculate the errors of numerical methods from the following relations:
\begin{equation}
\xi _F^{EX} = \mathop {\max }\limits_j \left( {\left| {x\left(
		{{t_j}} \right) - x_j^F} \right|} \right),\xi _{PC}^{EX} = \mathop {\max}\limits_j \left( {\left| {x\left( {{t_j}} \right) - x_j^{PC}} \right|} \right),
\label{Parovikeq23}
\end{equation}
where $ \xi _F^{EX}, \xi _ {PC} ^ {EX} $ ~ --- errors of ENFDS and ABM methods, orders of computational accuracy ~ --- $ p_F $ and $ p_ {PC} $, which calculated by the formula (\ref {Parovikeq20}).

\begin{table}[h!]
	\caption{Errors analysis of numerical methods for $ \beta = 1.8, \gamma = 0.9 $, taking into account the exact solution (\ref{Parovikeq22})}
	\centering
	\begin{ruledtabular}
	\begin{tabular}{|p{21pt}|p{28pt}|p{56pt}|p{64pt}|p{78pt}|p{71pt}|}
\parbox{21pt}{\centering $N$ } & \parbox{28pt}{\centering $\tau$} & \parbox{56pt}{\centering $\xi_F^{EX}$ } & \parbox{64pt}{\centering $\xi_{PC}^{EX}$} & \parbox{78pt}{\centering $p_F$} & \parbox{71pt}{\centering $p_{PC}$}\\
\hline
\parbox{21pt}{\centering 
			10
		} & \parbox{28pt}{\centering 
			1/10
		} & \parbox{56pt}{\centering 
			0.0561528
		} & \parbox{64pt}{\centering 
			0.0039560
		} & \parbox{78pt}{\centering 
			-
		} & \parbox{71pt}{\centering 
			-
		} \\
		\parbox{21pt}{\centering 
			20
		} & \parbox{28pt}{\centering 
			1/20
		} & \parbox{56pt}{\centering 
			0.0340440
		} & \parbox{64pt}{\centering 
			0.0009317
		} & \parbox{78pt}{\centering 
			0.721956135
		} & \parbox{71pt}{\centering 
			2.0861050
		} \\
		\parbox{21pt}{\centering 
			40
		} & \parbox{28pt}{\centering 
			1/40
		} & \parbox{56pt}{\centering 
			0.0195531
		} & \parbox{64pt}{\centering 
			0.0002362
		} & \parbox{78pt}{\centering 
			0.800003094
		} & \parbox{71pt}{\centering 
			1.9798565
		} \\
		\parbox{21pt}{\centering 
			80
		} & \parbox{28pt}{\centering 
			1/80
		} & \parbox{56pt}{\centering 
			0.0108341
		} & \parbox{64pt}{\centering 
			0.0000628
		} & \parbox{78pt}{\centering 
			0.851815427
		} & \parbox{71pt}{\centering 
			1.9100427
		} \\
		\parbox{21pt}{\centering 
			160
		} & \parbox{28pt}{\centering 
			1/160
		} & \parbox{56pt}{\centering 
			0.0058603
		} & \parbox{64pt}{\centering 
			0.0000171
		} & \parbox{78pt}{\centering 
			0.886538801
		} & \parbox{71pt}{\centering 
			1.8748556
		} \\
		\parbox{21pt}{\centering 
			320
		} & \parbox{28pt}{\centering 
			1/320
		} & \parbox{56pt}{\centering 
			0.0031176
		} & \parbox{64pt}{\centering 
			0.0000047
		} & \parbox{78pt}{\centering 
			0.910518254
		} & \parbox{71pt}{\centering 
			1.8577979
		} \\
		\parbox{21pt}{\centering 
			640
		} & \parbox{28pt}{\centering 
			1/640
		} & \parbox{56pt}{\centering 
			0.0016380
		} & \parbox{64pt}{\centering 
			0.0000014
		} & \parbox{78pt}{\centering 
			0.92851899
		} & \parbox{71pt}{\centering 
			1.7975624
		} \\
	\end{tabular}
\end{ruledtabular}	
\end{table} 

From Table 1 shows that with an increase in the nodes of the computational grid, the errors of numerical methods decrease. It can be seen that the ABM method converges faster than the JANCRS method. Based on (\ref {Parovikeq9}), the ENFDS approximation is of the first order. According to the relation (\ref{Parovikeq18}) for Example 1, the order of the ABM method.

From Table 1, we also see that with an increase in the nodes of the computational grid, the computational accuracy increases and tends to unity, and $ p_ {PC} $ ~ --- almost immediately reached the value 1.9. On the one hand, this confirms the validity of Lemma 3 and Theorem 3, and on the other hand, the faster convergence and high accuracy of the ABM method in comparison with the ENFDS method.

As a comparison with the results from Table 1, we calculate the errors and orders of magnitude of the computational accuracy of numerical methods taking into account the Runge rule (Table 2).

\begin{table}[h!]
	\caption{Errors analysis of numerical methods for $\beta = 1.8, \gamma = 0.9 $ using the double recalculation method (\ref{Parovikeq19})}
	\centering
\begin{ruledtabular}
	\begin{tabular}{|p{21pt}|p{28pt}|p{49pt}|p{64pt}|p{78pt}|p{71pt}|}      \hline
		\parbox{21pt}{\centering $N$ } & \parbox{28pt}{\centering $\tau$} & \parbox{56pt}{\centering $\xi_F^{EX}$ } & \parbox{64pt}{\centering $\xi_{PC}^{EX}$} & \parbox{78pt}{\centering $p_F$} & \parbox{71pt}{\centering $p_{PC}$} \\
		\hline
		\parbox{21pt}{\centering 
			10
		} & \parbox{28pt}{\centering 
			1/10
		} & \parbox{49pt}{\centering 
			0.025634
		} & \parbox{64pt}{\centering 
			0.0011076
		} & \parbox{78pt}{\centering 
			\textbf{-}
		} & \parbox{71pt}{\centering 
			\textbf{-}
		} \\
		\parbox{21pt}{\centering 
			20
		} & \parbox{28pt}{\centering 
			1/20
		} & \parbox{49pt}{\centering 
			0.015479
		} & \parbox{64pt}{\centering 
			0.0002546
		} & \parbox{78pt}{\centering 
			0.727709917
		} & \parbox{71pt}{\centering 
			2.120946173
		} \\
		\parbox{21pt}{\centering 
			40
		} & \parbox{28pt}{\centering 
			1/40
		} & \parbox{49pt}{\centering 
			0.008987
		} & \parbox{64pt}{\centering 
			0.0000635
		} & \parbox{78pt}{\centering 
			0.784362481
		} & \parbox{71pt}{\centering 
			2.004652815
		} \\
		\parbox{21pt}{\centering 
			80
		} & \parbox{28pt}{\centering 
			1/80
		} & \parbox{49pt}{\centering 
			0.005045
		} & \parbox{64pt}{\centering 
			0.0000167
		} & \parbox{78pt}{\centering 
			0.832949245
		} & \parbox{71pt}{\centering 
			1.923159876
		} \\
		\parbox{21pt}{\centering 
			160
		} & \parbox{28pt}{\centering 
			1/160
		} & \parbox{49pt}{\centering 
			0.002761
		} & \parbox{64pt}{\centering 
			0.0000045
		} & \parbox{78pt}{\centering 
			0.86953863
		} & \parbox{71pt}{\centering 
			1.881091919
		} \\
		\parbox{21pt}{\centering 
			320
		} & \parbox{28pt}{\centering 
			1/320
		} & \parbox{49pt}{\centering 
			0.001484
		} & \parbox{64pt}{\centering 
			0.0000013
		} & \parbox{78pt}{\centering 
			0.895926947
		} & \parbox{71pt}{\centering 
			1.860912083
		} \\
		\parbox{21pt}{\centering 
			640
		} & \parbox{28pt}{\centering 
			1/640
		} & \parbox{49pt}{\centering 
			0.000786
		} & \parbox{64pt}{\centering 
			0.0000004
		} & \parbox{78pt}{\centering 
			0.916645745
		} & \parbox{71pt}{\centering 
			1.78891439
		} \\
	\end{tabular}
\end{ruledtabular}	
\end{table}

From Table 2 that the double recalculation method based on the Runge rule is in good agreement with the results from Table 1.

In the classical case, when the ABM method has the second order, and the ENFDS is still the same first order of approximation. Numerical analysis of errors for Example 1 is given in Table 3.

\begin{table}[h!]
	\centering	
	\caption{Errors analysis of numerical methods for $ \beta = 2, \gamma = 1 $ by the double recalculation method (\ref{Parovikeq19})}	
\begin{ruledtabular}
	\begin{tabular}{|p{21pt}|p{35pt}|p{64pt}|p{78pt}|p{78pt}|p{71pt}|}
		\hline
		\parbox{21pt}{\centering $N$ } & \parbox{28pt}{\centering $\tau$} & \parbox{56pt}{\centering $\xi_F$} & \parbox{64pt}{\centering $\xi_{PC}$} & \parbox{78pt}{\centering $p_F$} & \parbox{71pt}{\centering $p_{PC}$}\\
		\hline
		\parbox{21pt}{\centering 
			10
		} & \parbox{35pt}{\centering 
			1/10
		} & \parbox{64pt}{\centering 
			0.00758
		} & \parbox{78pt}{\centering 
			0.001049
		} & \parbox{49pt}{\centering 
			\textbf{-}
		} & \parbox{49pt}{\centering 
			\textbf{-}
		} \\
		\parbox{21pt}{\centering 
			20
		} & \parbox{35pt}{\centering 
			1/20
		} & \parbox{64pt}{\centering 
			0.00253
		} & \parbox{78pt}{\centering 
			0.000232
		} & \parbox{49pt}{\centering 
			1.585196
		} & \parbox{49pt}{\centering 
			2.179457
		} \\
		\parbox{21pt}{\centering 
			40
		} & \parbox{35pt}{\centering 
			1/40
		} & \parbox{64pt}{\centering 
			0.00093
		} & \parbox{78pt}{\centering 
			0.000054
		} & \parbox{49pt}{\centering 
			1.445025
		} & \parbox{49pt}{\centering 
			2.093444
		} \\
		\parbox{21pt}{\centering 
			80
		} & \parbox{35pt}{\centering 
			1/80
		} & \parbox{64pt}{\centering 
			0.00038
		} & \parbox{78pt}{\centering 
			0.000013
		} & \parbox{49pt}{\centering 
			1.296494
		} & \parbox{49pt}{\centering 
			2.046768
		} \\
		\parbox{21pt}{\centering 
			160
		} & \parbox{35pt}{\centering 
			1/160
		} & \parbox{64pt}{\centering 
			0.00017
		} & \parbox{78pt}{\centering 
			0.000003
		} & \parbox{49pt}{\centering 
			1.177333
		} & \parbox{49pt}{\centering 
			2.023954
		} \\
		\parbox{21pt}{\centering 
			320
		} & \parbox{35pt}{\centering 
			1/320
		} & \parbox{64pt}{\centering 
			0.00008
		} & \parbox{78pt}{\centering 
			0.0000008007
		} & \parbox{49pt}{\centering 
			1.100742
		} & \parbox{49pt}{\centering 
			2.012096
		} \\
		\parbox{21pt}{\centering 
			640
		} & \parbox{35pt}{\centering 
			1/640
		} & \parbox{64pt}{\centering 
			0.00003739
		} & \parbox{78pt}{\centering 
			0.0000001994
		} & \parbox{49pt}{\centering 
			1.058277
		} & \parbox{49pt}{\centering 
			2.005656
		} \\
	\end{tabular}
\end{ruledtabular}
\end{table}

From Table 3 that the orders of computational accuracy and tend to the orders of approximation 1 and 2 for the corresponding numerical methods.

\begin{figure}[h!]
	\centering
	\includegraphics[scale=0.65]{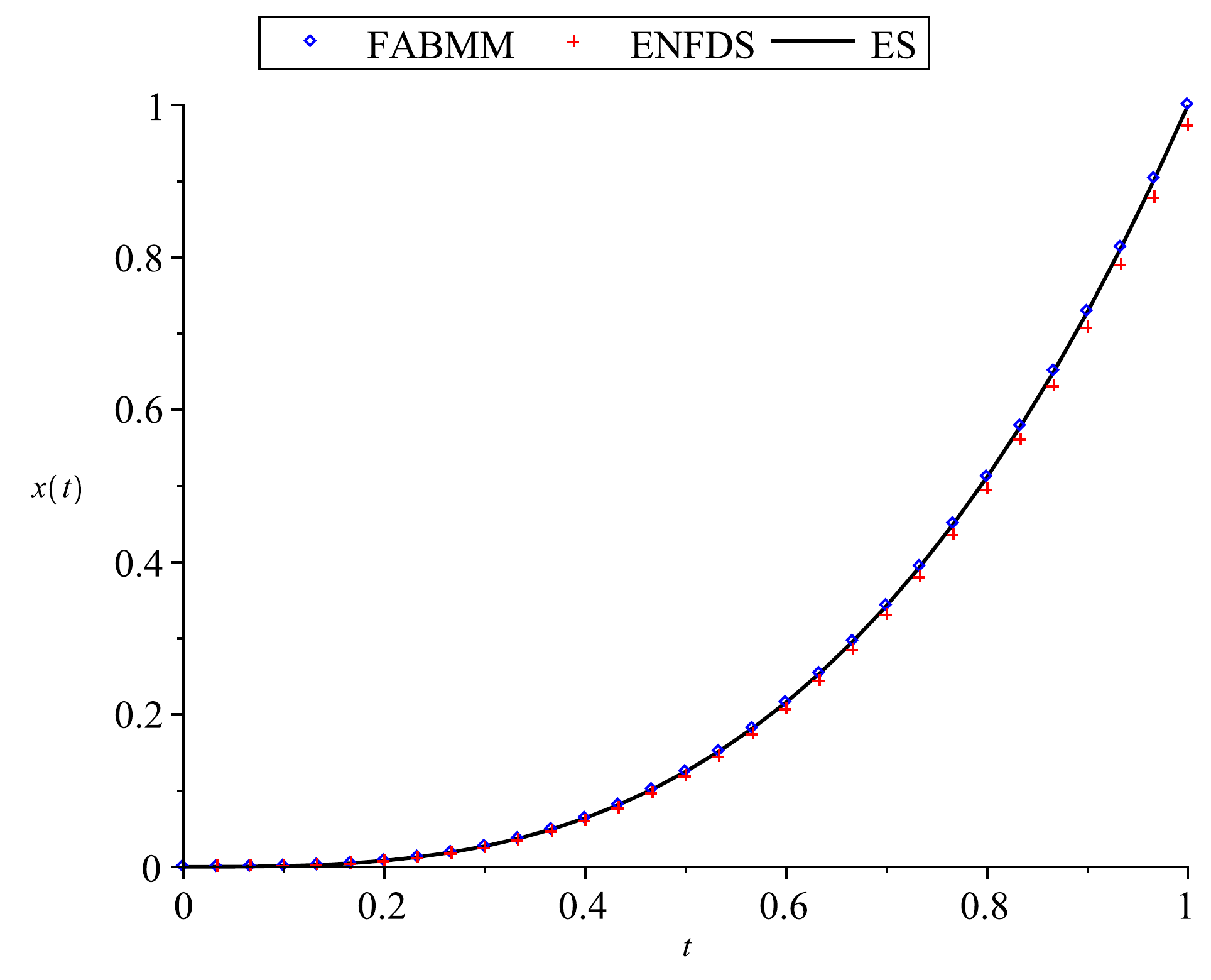}
	\caption{Numerical results using the Adams-Bashfort-Moulton (ABM) and Explicit Nonlocal Finite Difference Scheme (ENFDS) methods compared to the exact solution (ES) for Example 1 at $ N = 20 $}\label{Parovikfig1}
\end{figure}

Figure \ref{Parovikfig1} illustrates the effectiveness of the ABM method in comparison with the ENFDS. It can be seen that, due to the fast convergence and higher accuracy, the ABM method gives the best result. 

However, the ABM method does not always work better than explicit methods like ENFDS. Let us show this with the following example.

\textbf{Example 2.} Consider an analogue of the harmonic oscillator (AHO). To do this, in the Cauchy problem (\ref {Parovikeq3}) choose:
	\[
	\omega \left( t \right) = \omega _0^\beta ,f\left( t \right) = \frac{{\mu
			{t^\delta }}}{{\Gamma \left( {1 + \delta } \right)}},\lambda  = 0.
	\]
Then we come to the following problem:
\begin{equation}
\partial _{0t}^\beta x\left( \tau  \right) + \omega _0^\beta x\left( t \right) =
\frac{{\mu {t^\delta }}}{{\Gamma \left( {1 + \delta } \right)}},x\left( 0 \right)
= {\alpha _1},\dot x\left( 0 \right) = {\alpha _2}.
\label{Parovikeq24}
\end{equation}

The solution to the Cauchy problem (\ref{Parovikeq24}) is the function:
\begin{equation}
x\left( t \right) = {\alpha _1}{E_{\beta ,1}}\left( { - {{\left( {{\omega_0}t} \right)}^\beta }} \right) + {\alpha _2}t{E_{\beta, 2}}\left( { - {{\left({{\omega _0}t} \right)}^\beta }} \right) + \mu {t^{\beta  + g}}{E_{\beta,\,\beta + \delta  + 1}}\left( { - {{\left( {{\omega _0}t} \right)}^\beta }} \right), \delta  >  - \beta,
\label{Parovikeq25}   
\end{equation}	
where ${E_{\alpha ,\beta }}\left( z \right) = \sum\limits_{k = 0}^\infty {\dfrac{{{z^k}}}{{\Gamma \left( {\alpha k + \beta } \right)}}} $~--- two-parameter Mittag-Leffler function \cite{Kilbas_2006}.

Let's consider the numerical methods ABM and ENFDS for solving the Cauchy problem (\ref{Parovikeq24}) and compare them with the exact solution (\ref{Parovikeq25}). First, consider the usual harmonic oscillator ($ \beta = 2 $).

The parameter values in the problem (\ref{Parovikeq25}) are chosen as follows: $\delta  = 0.3,\mu  =
0.5,{\omega _0} = 2,{\alpha _1} = {\alpha _2} = 1,t \in \left[ {0,1} \right]$.

\begin{figure}[h!]
	\centering	
	\includegraphics[scale=0.65]{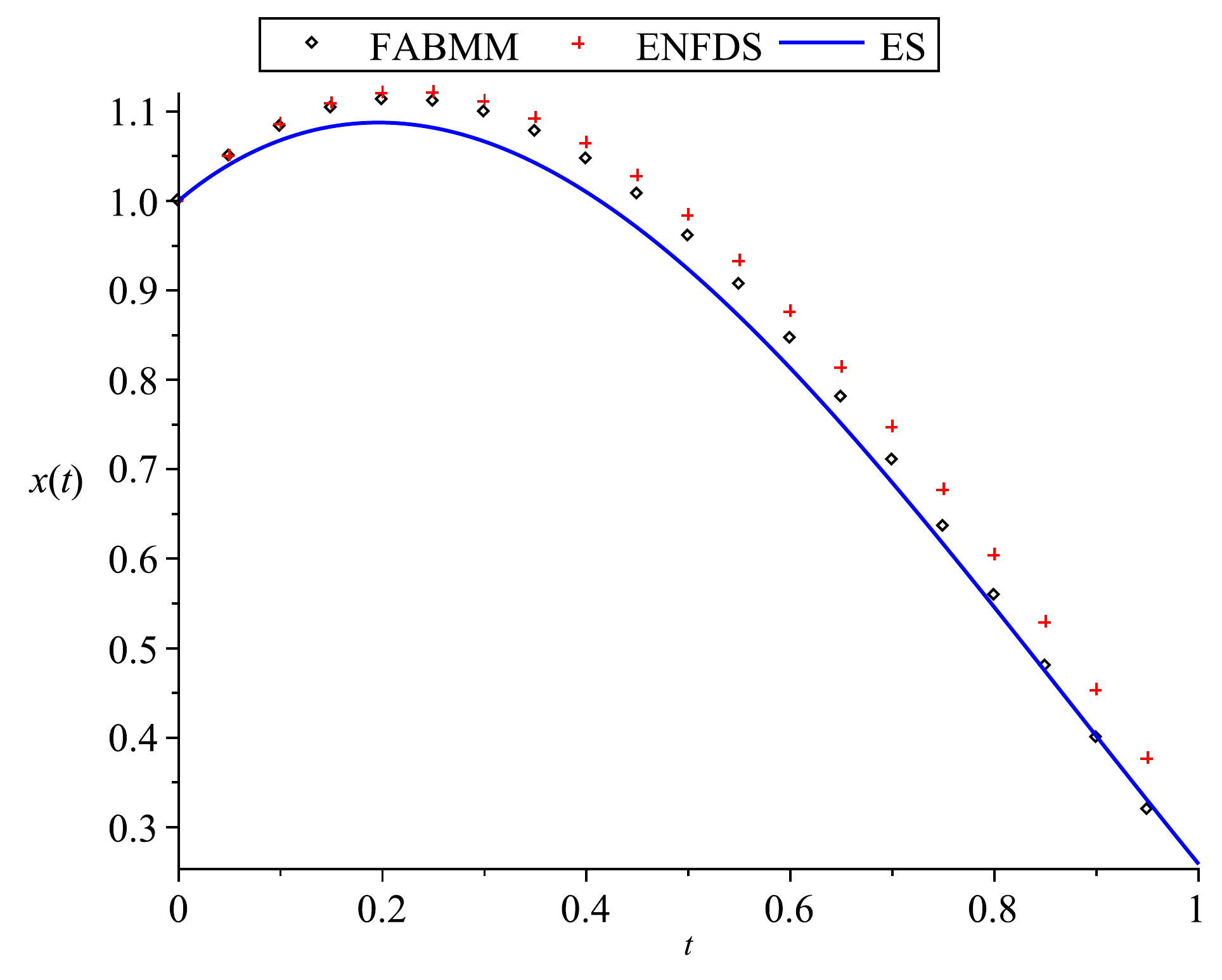}
	\caption{Numerical results using Adams-Bashfort-Moulton (ABM) and Explicit Nonlocal Finite Difference Scheme (ENFDS) versus Exact Solution (ES) for Example 2 at $ N = 20 $}
\end{figure}

We see in Figure 2 that the ABM method works in much the same way as the ENFDS method. The order of computational accuracy is shown in the following Table 4. In the case of AHO, the error analysis is given in Table 5.

\begin{table}[h!]
	\centering	
	\caption{Errors analysis of numerical methods for Example 2 ($\beta=2$). }
\begin{ruledtabular}
	\begin{tabular}{|p{21pt}|p{35pt}|p{77pt}|p{78pt}|p{70pt}|p{70pt}|}
		\hline
		\parbox{21pt}{\centering $N$ } & \parbox{28pt}{\centering $\tau$} & \parbox{56pt}{\centering $\xi_F^{EX}$ } & \parbox{64pt}{\centering $\xi_{PC}^{EX}$} & \parbox{78pt}{\centering $p_F$} & \parbox{71pt}{\centering $p_{PC}$ } \\
		\hline 
		\parbox{21pt}{\centering 
			10
		} & \parbox{35pt}{\centering 
			1/10
		} & \parbox{77pt}{\centering 
			0.0962543446
		} & \parbox{78pt}{\centering 
			0.0974594278
		} & \parbox{70pt}{\centering 
			\textbf{-}
		} & \parbox{70pt}{\centering 
			\textbf{-}
		} \\
		\parbox{21pt}{\centering 
			20
		} & \parbox{35pt}{\centering 
			1/20
		} & \parbox{77pt}{\centering 
			0.0483529972
		} & \parbox{78pt}{\centering 
			0.0450644425
		} & \parbox{70pt}{\centering 
			0.993246342
		} & \parbox{70pt}{\centering 
			1.112812209
		} \\
		\parbox{21pt}{\centering 
			40
		} & \parbox{35pt}{\centering 
			1/40
		} & \parbox{77pt}{\centering 
			0.0243287977
		} & \parbox{78pt}{\centering 
			0.0216596566
		} & \parbox{70pt}{\centering 
			0.990940293
		} & \parbox{70pt}{\centering 
			1.056979173
		} \\
		\parbox{21pt}{\centering 
			80
		} & \parbox{35pt}{\centering 
			1/80
		} & \parbox{77pt}{\centering 
			0.0122257234
		} & \parbox{78pt}{\centering 
			0.0106471450
		} & \parbox{70pt}{\centering 
			0.992745193
		} & \parbox{70pt}{\centering 
			1.024543742
		} \\
		\parbox{21pt}{\centering 
			160
		} & \parbox{35pt}{\centering 
			1/160
		} & \parbox{77pt}{\centering 
			0.0061384979
		} & \parbox{78pt}{\centering 
			0.0052867956
		} & \parbox{70pt}{\centering 
			0.993962258
		} & \parbox{70pt}{\centering 
			1.010001173
		} \\
		\parbox{21pt}{\centering 
			320
		} & \parbox{35pt}{\centering 
			1/320
		} & \parbox{77pt}{\centering 
			0.0030783418
		} & \parbox{78pt}{\centering 
			0.0026378633
		} & \parbox{70pt}{\centering 
			0.995732241
		} & \parbox{70pt}{\centering 
			1.003023747
		} \\
		\parbox{21pt}{\centering 
			640
		} & \parbox{35pt}{\centering 
			1/640
		} & \parbox{77pt}{\centering 
			0.0015461838
		} & \parbox{78pt}{\centering 
			0.0013190516
		} & \parbox{70pt}{\centering 
			0.993441601
		} & \parbox{70pt}{\centering 
			0.9998688
		} \\
	\end{tabular}
\end{ruledtabular}	
\end{table}

\begin{table}[h!]
	\centering	
	\caption{Errors analysis of numerical methods for Example 2 ($\beta=1.8$)}
\begin{ruledtabular}	
	\begin{tabular}{|p{21pt}|p{35pt}|p{77pt}|p{78pt}|p{70pt}|p{70pt}|}
		\hline
		\parbox{21pt}{\centering $N$ } & \parbox{28pt}{\centering $\tau$} & \parbox{56pt}{\centering $\xi_F^{EX}$ } & \parbox{64pt}{\centering $\xi_{PC}^{EX}$} & \parbox{78pt}{\centering $p_F$} & \parbox{71pt}{\centering  $p_{PC}$}\\
		\hline
		\parbox{21pt}{\centering 
			10
		} & \parbox{35pt}{\centering 
			1/10
		} & \parbox{77pt}{\centering 
			0.1220802888
		} & \parbox{78pt}{\centering 
			0.0833919341
		} & \parbox{70pt}{\centering 
			\textbf{-}
		} & \parbox{70pt}{\centering 
			\textbf{-}
		} \\
		\parbox{21pt}{\centering 
			20
		} & \parbox{35pt}{\centering 
			1/20
		} & \parbox{77pt}{\centering 
			0.0632933442
		} & \parbox{78pt}{\centering 
			0.0376498054
		} & \parbox{70pt}{\centering 
			0.947704578
		} & \parbox{70pt}{\centering 
			1.147265441
		} \\
		\parbox{21pt}{\centering 
			40
		} & \parbox{35pt}{\centering 
			1/40
		} & \parbox{77pt}{\centering 
			0.0326304379
		} & \parbox{78pt}{\centering 
			0.0179210448
		} & \parbox{70pt}{\centering 
			0.955835448
		} & \parbox{70pt}{\centering 
			1.070987659
		} \\
		\parbox{21pt}{\centering 
			80
		} & \parbox{35pt}{\centering 
			1/80
		} & \parbox{77pt}{\centering 
			0.0167613673
		} & \parbox{78pt}{\centering 
			0.0087578726
		} & \parbox{70pt}{\centering 
			0.961078508
		} & \parbox{70pt}{\centering 
			1.033002381
		} \\
		\parbox{21pt}{\centering 
			160
		} & \parbox{35pt}{\centering 
			1/160
		} & \parbox{77pt}{\centering 
			0.0085801061
		} & \parbox{78pt}{\centering 
			0.0043373867
		} & \parbox{70pt}{\centering 
			0.966072448
		} & \parbox{70pt}{\centering 
			1.013754391
		} \\
		\parbox{21pt}{\centering 
			320
		} & \parbox{35pt}{\centering 
			1/320
		} & \parbox{77pt}{\centering 
			0.0043785606
		} & \parbox{78pt}{\centering 
			0.0021615271
		} & \parbox{70pt}{\centering 
			0.970538809
		} & \parbox{70pt}{\centering 
			1.004775148
		} \\
		\parbox{21pt}{\centering 
			640
		} & \parbox{35pt}{\centering 
			1/640
		} & \parbox{77pt}{\centering 
			0.0022275871
		} & \parbox{78pt}{\centering 
			0.0010803221
		} & \parbox{70pt}{\centering 
			0.974974836
		} & \parbox{70pt}{\centering 
			1.000589405
		} \\
	\end{tabular}
\end{ruledtabular}
\end{table}

We see that the order of the computational accuracy of the methods tends to unity, i.e. the accuracy of the ABM method is similar to that of the ENFDS method, although the convergence rate is higher.

\section{Conclusion}

In this work, a numerical analysis of the Cauchy problem for a model class of fractional linear oscillators was carried out. The methods of ENFDS and ABM were considered as numerical methods. It has been shown on test examples that the ABM method converges faster and is more accurate (not always) than the ENFDS method.

The study of FLO with variable memory is also of interest. In this case, derivatives of fractional variable orders appear in the model equation of the Cauchy problem (\ref{Parovikeq3}). Some aspects of the numerical analysis of such a generalized Cauchy problem (\ref{Parovikeq3}) were considered by the author in the work \cite{Parovik_2016}.

\nocite{*}
\bibliography{Parovik}

\end{document}